\documentclass[11pt]{article}
\usepackage{amsfonts}

\usepackage[dvips]{graphicx}
\usepackage{amsmath}

\newtheorem{theorem}{Theorem}[subsection]

\newtheorem{example}[theorem]{Example}

\newtheorem{lemma}[theorem]{Lemma}

\newtheorem{proposition}[theorem]{Proposition}

\newenvironment{proof}[1][Proof]{\textbf{#1.} }{\ \rule{0.5em}{0.5em}}

\input{tcilatex}

\begin{document}

\title{Stabilized plethysms for the classical Lie groups }
\author{C\'{e}dric Lecouvey \\
%EndAName
Laboratoire de Math\'{e}matiques Pures et Appliqu\'{e}es Joseph Liouville\\
B.P. 699 62228 Calais Cedex\\
Cedric.Lecouvey@lmpa.univ-littoral.fr}
\date{}
\maketitle

\begin{abstract}
The plethysms of the Weyl characters associated to a classical Lie group by
the symmetric functions stabilize in large rank. In the case of a power sum
plethysm, we prove that the coefficients of the decomposition of this
stabilized form on the basis of Weyl characters are branching coefficients
which can be determined by a simple algorithm.\ This generalizes in
particular some classical results by Littlewood on the power sum plethysms
of Schur functions. We also establish explicit formulas for the outer
multiplicities appearing in the decomposition of the tensor square of any
irreducible finite dimensional module into its symmetric and antisymmetric
parts. These multiplicities can notably be expressed in terms of the
Littlewood-Richardson coefficients.
\end{abstract}

\section{Introduction}

This paper is concerned with the plethysms of the Weyl characters associated
to classical Lie groups by the symmetric functions. Let $\frak{g}$ be a
classical Lie group with rank $n,$ and $\lambda $ a partition.\ We denote by 
$s_{\lambda }^{\frak{g}}$ the Weyl character of the $\frak{g}$-module $V^{%
\frak{g}}(\lambda )$ (see Section \ref{Sec_Back}). Consider $f$ a symmetric
function of degree $d$ and suppose $n>dl(\lambda )$ where $l(\lambda )$ is
the number of non-zero parts of $\lambda $ It follows from results by
Littlewood \cite{Li2} that the plethysm $f\circ s_{\lambda }^{\frak{g}}$ of
the Weyl character $s_{\lambda }^{\frak{g}}$ by $f$ decomposes on the basis $%
\{s_{\mu }^{\frak{g}}\mid \mu \in \mathcal{P}_{n}\}$ with coefficients which
do not depend on $n$.\ When $f=p_{\ell }$ is the power sum of degree $\ell ,$
we establish that the coefficients so obtained are branching coefficients
corresponding to the restriction to certain Levi subgroups (Theorem \ref
{th_a=bL}). Suppose $n>\ell l(\lambda )$ and set 
\begin{equation*}
p_{\ell }\circ s_{\lambda }^{\frak{g}}=\sum_{\mu }a_{\lambda ,\mu }^{\frak{%
g,\ell }}s_{\mu }^{\frak{g}}.
\end{equation*}
For $\frak{g=gl}_{n},$ it is well known, by an algorithm due to Littlewood 
\cite{Li0}, that the coefficients $a_{\lambda ,\mu }^{\frak{gl}_{n}\frak{%
,\ell }}$ can, up to a sign, be expressed as a sum of products of
Littlewood-Richardson coefficients.\ They are then obtained from the $\ell $%
-quotient of the partition $\mu .\;$We give a similar algorithm for
computing the coefficients $a_{\lambda ,\mu }^{\frak{g,\ell }}$ when $\frak{%
g=so}_{2n+1},\frak{sp}_{2n}$ or $\frak{so}_{2n}.$ This algorithm was
originally introduced in \cite{lec} to decompose the plethysms $p_{\ell
}\circ s_{\lambda }^{\frak{so}_{2n+1}}$ on the basis of Weyl characters for
any integers $n\geq 2$ and $\ell \geq 1$ (that is, with no restrictive
conditions on the rank $n$). Although similar procedures also exist for $%
\frak{g=sp}_{2n}$ or $\frak{so}_{2n}$ when $\ell $ is odd, our method failed
for the even power sum plethysms on the Weyl characters of type $C_{n}$ or $%
D_{n}$. In the present paper, we show that this difficulty can be overcome
by considering stabilized power sum plethysms, i.e.\ by assuming that $%
n>\ell l(\lambda ).$ Under this hypothesis, one has indeed $a_{\lambda ,\mu
}^{\frak{so}_{2n+1\frak{,\ell }}}=a_{\lambda ,\mu }^{\frak{so}_{2n\frak{%
,\ell }}}$ and $a_{\lambda ,\mu }^{\frak{sp}_{2n\frak{,\ell }}}=(-1)^{\left|
\lambda \right| }(-1)^{\ell -1}a_{\lambda ^{\prime },\mu ^{\prime }}^{\frak{%
so}_{2n+1}\frak{,\ell }}$.\ So it suffices to consider the coefficients $%
a_{\lambda ,\mu }^{\frak{so}_{2n+1}\frak{,\ell }}$ for which there exists an
algorithm in both cases $\ell $ even and $\ell $ odd.

\noindent In Proposition \ref{prop-symanti}, we use our expression of the
coefficients $a_{\lambda ,\mu }^{\frak{g,}2}$ as branching coefficients, to
derive explicit formulas giving the decompositions of the symmetric and
antisymmetric parts of $V^{\frak{g}}(\lambda )^{\otimes 2}$ in their
irreducible components when $n>2l(\lambda ).$ The corresponding
multiplicities can then be expressed in terms of the Littlewood-Richardson
coefficients and give an alternative to analogous formulas introduced
without a complete proof by Littlewood in \cite{Li2}.

\bigskip

\noindent The paper is organized as follows.\ In Section $2,$ we recall some
basics on the representation theory of the classical Lie groups. Section $3$
is concerned with plethysms $f\circ s_{\lambda }^{\frak{g}}$ and their
stabilization in large rank.\ Most of the material of this section can be
found in \cite{Li0}, \cite{Li}, \cite{Li2} and \cite{mac}.\ In Section $4$,
we describe the algorithm of \cite{lec} which permits to compute the
plethysms $p_{\ell }\circ s_{\lambda }^{\frak{so}_{2n+1}}$ for any positive
integer $\ell .$ We then state Theorem \ref{th_a=bL}.\ Finally, in Section $%
5,$ we express the multiplicities $a_{\lambda ,\mu }^{\frak{g,}2}$ in terms
of the Littlewood-Richardson coefficients.

\bigskip

\noindent \textbf{Acknowledgments: }The author wants to thank the anonymous
referees for having pointed out some mistakes and inacurrencies in a
previous version of this paper.\ In particular, the stabilization phenomenon
explained in Section \ref{subsec_rq} emerge now naturally from the
algorithms of Section \ref{subsec-algoeven} and \ref{subsec-algoodd}. This
yields to stabilization conditions stated in Theorem \ref{Th_LI}.

\section{Background on classical Lie groups\label{Sec_Back}}

\subsection{Root systems and Weyl groups}

In the sequel $G$ is one of the complex Lie groups $Sp_{2n},SO_{2n+1}$ or $%
SO_{2n}$ and $\frak{g}$ is its Lie algebra.\ We follow the convention of 
\cite{KT} to realize $G$ as a subgroup of $GL_{N}$ and $\frak{g}$\ as a
subalgebra of $\frak{gl}_{N}$ where 
\begin{equation*}
N=\left\{ 
\begin{tabular}{l}
$n$ when $G=GL_{n}$ \\ 
$2n$ when $G=Sp_{2n}$ \\ 
$2n+1$ when $G=SO_{2n+1}$ \\ 
$2n$ when $G=SO_{2n}.$%
\end{tabular}
\right. .
\end{equation*}
\noindent Let $d_{N}$ be the linear subspace of $\frak{gl}_{N}$ consisting
of the diagonal matrices.\ For any $i\in I_{n}=\{1,...,n\},$ write $%
\varepsilon _{i}$ for the linear map $\varepsilon _{i}:d_{N}\rightarrow 
\mathbb{C}$ such that $\varepsilon _{i}(D)=\delta _{n-i+1}$ for any diagonal
matrix $D$ whose $(i,i)$-coefficient is $\delta _{i}.$ Then $(\varepsilon
_{1},...,\varepsilon _{n})$ is an orthonormal basis of the Euclidean space $%
\frak{h}_{\mathbb{R}}^{\ast }$ (the real part of $\frak{h}^{\ast }).$ Let $%
(\cdot ,\cdot )$ be the corresponding nondegenerate symmetric bilinear form
defined on $\frak{h}_{\mathbb{R}}^{\ast }$.\ Write $R$ for the root system
associated to $G.\;$For any $\alpha \in R$ we set $\alpha ^{\vee }=\frac{%
\alpha }{(\alpha ,\alpha )}$. The Lie algebra $\frak{g}$ admits the diagonal
decomposition $\frak{g}=\frak{h}\oplus \tcoprod_{\alpha \in R}\frak{g}%
_{\alpha }.\;$We take for the set of positive roots: 
\begin{equation*}
\left\{ 
\begin{tabular}{l}
$R^{+}=\{\varepsilon _{j}-\varepsilon _{i}\text{ with }1\leq i<j\leq n\}%
\text{ for the root system }A_{n-1}$ \\ 
$R^{+}=\{\varepsilon _{j}-\varepsilon _{i},\varepsilon _{j}+\varepsilon _{i}%
\text{ with }1\leq i<j\leq n\}\cup \{\varepsilon _{i}\text{ with }1\leq
i\leq n\}\text{ for the root system }B_{n}$ \\ 
$R^{+}=\{\varepsilon _{j}-\varepsilon _{i},\varepsilon _{j}+\varepsilon _{i}%
\text{ with }1\leq i<j\leq n\}\cup \{2\varepsilon _{i}\text{ with }1\leq
i\leq n\}\text{ for the root system }C_{n}$ \\ 
$R^{+}=\{\varepsilon _{j}-\varepsilon _{i},\varepsilon _{j}+\varepsilon _{i}%
\text{ with }1\leq i<j\leq n\}\text{ for the root system }D_{n}$%
\end{tabular}
\right. .
\end{equation*}
For any $i\in I_{n}$, we write $\overline{i}$ for $-i$.\ The Weyl group $W$
of the Lie group $G$ is the subgroup of the permutation group of the set $%
J_{n}=\{\overline{n},...,\overline{2},\overline{1},1,2,...,n\}$\ generated
by the permutations 
\begin{equation*}
\left\{ 
\begin{tabular}{l}
$s_{i}=(i,i+1)(\overline{i},\overline{i+1}),$ $i=1,...,n-1$ and $s_{0}=(1,%
\overline{1})$ $\text{for the root systems }B_{n}$ and $C_{n}$ \\ 
$s_{i}=(i,i+1)(\overline{i},\overline{i+1}),$ $i=1,...,n-1$ and $%
s_{0}^{\prime }=(1,\overline{2})(2,\overline{1})$ $\text{for the root system 
}D_{n}$%
\end{tabular}
\right.
\end{equation*}
where for $a\neq b$, $(a,b)$ is the simple transposition which switches $a$
and $b.$ We identify the subgroup of $W$ generated by $s_{i}=(i,i+1)(%
\overline{i},\overline{i+1}),$ $i=1,...,n-1$ with the symmetric group $%
S_{n}. $ We denote by $l$ the length function corresponding to the above set
of generators. For any $w\in W,$ we set $\varepsilon (w)=(-1)^{l(w)}.\;$The
action of $w\in W$ on $\beta =(\beta _{1},...,\beta _{n})\in \frak{h}_{%
\mathbb{R}}^{\ast }$ is defined by 
\begin{equation*}
w\cdot (\beta _{1},...,\beta _{n})=(\beta _{1}^{w^{-1}},...,\beta
_{n}^{w^{-1}})
\end{equation*}
where $\beta _{i}^{w}=\beta _{w(i)}$ if $w(i)\in \{1,...,n\}$ and $\beta
_{i}^{w}=-\beta _{w(\overline{i})}$ otherwise. We denote by $\rho $ the half
sum of the positive roots of $R^{+}$. For any $x\in J_{n},$ we set $%
\overline{\overline{x}}=x$ and $\left| x\right| =x$ if $x$ is unbarred, $%
\left| x\right| =\overline{x}$ otherwise.

\bigskip

\noindent A partition of length $m$ is a weakly \textit{increasing} sequence
of $m$ nonnegative integers.\ Denote by $\mathcal{P}_{m}$ the set of
partitions with at most $m$ parts. Given $\lambda \in \mathcal{P}_{m}$, $%
\lambda ^{\prime }$ is its conjugate partition and $l(\lambda )$ the number
of nonzero parts in $\lambda $.\ Set $\mathcal{P}=\cup _{m\geq 0}\mathcal{P}%
_{m}.\;$For $G=Sp_{2n}$ or $SO_{2n+1}$ and $\lambda \in \mathcal{P}_{n},$
denote by $V^{\frak{g}}(\lambda )$ the irreducible finite dimensional
representation of $G$ of highest weight $\lambda .\;$For $G=SO_{2n}$, we
define $V^{\frak{so}_{2n}}(\lambda )$ similarly when $\lambda _{1}=0$ and we
write $V^{\frak{so}_{2n}}(\lambda )$ for the direct sum of the two
irreducible representations of highest weights $\lambda =(\lambda
_{1},...,\lambda _{n-1},\lambda _{n})$ and $\overline{\lambda }=(-\lambda
_{1},...,\lambda _{n-1},\lambda _{n})$ when $\lambda _{n}\neq 0$. This means
that $V^{\frak{so}_{2n}}(\lambda )$ is in fact the irreducible
representation of $O_{2n}$ associated to the partition $\lambda $.

\noindent We shall also need the irreducible rational representations of $%
GL_{n}$.\ They are indexed by the $n$-tuples 
\begin{equation}
(\gamma ^{-},\gamma ^{+})=(-\gamma _{q}^{-},...,-\gamma _{1}^{-},\gamma
_{1}^{+},\gamma _{2}^{+},...,\gamma _{p}^{+})  \label{jamma+-}
\end{equation}
where $\gamma ^{+}=(\gamma _{1}^{+},\gamma _{2}^{+},...,\gamma _{p}^{+})$
and $\gamma ^{-}=(\gamma _{1}^{-},...,\gamma _{q}^{-})$ are partitions of
length $p$ and $q$ such that $p+q=n.$ Write $\widetilde{\mathcal{P}}_{n}$
for the set of such $n$-tuples and denote also by $V^{\frak{gl}_{n}}(\gamma
) $ the irreducible rational representation of $\frak{gl}_{n}$ of highest
weight $\gamma =(\gamma ^{-},\gamma ^{+})\in \widetilde{\mathcal{P}}_{n}.$
For any $\gamma =(\gamma ^{-},\gamma ^{+})\in \widetilde{\mathcal{P}}_{n},$
we set $\left| \gamma \right| =\sum \gamma _{i}^{-}+\sum \gamma _{i}^{+}.$

\bigskip

\noindent Write $s_{\lambda}^{\frak{gl}_{n}}$ for the Weyl character (Schur
function) of the finite-dimensional $\frak{gl}_{n}$-module $V^{\frak{gl}%
_{n}}(\lambda)$ of highest weight $\lambda.$ The character ring of $GL_{n}$
is $\Lambda_{n}=\mathbb{Z[}x_{1},...,x_{n}]^{\mathrm{sym}}$ the ring of
symmetric functions in $n$ variables.\ 

\noindent For any $\lambda \in \mathcal{P}_{n},$ we denote by $s_{\lambda }^{%
\frak{g}}$ the Weyl character of $V^{\frak{g}}(\lambda )$. Let $\mathcal{R}^{%
\frak{g}}$ be the $\mathbb{Z}$-algebra with basis $\{s_{\lambda }^{\frak{g}%
}\mid \lambda \in \mathcal{P}_{n}\}.$

\noindent Consider $P$ a parabolic subgroup of $G$ and $L$ its Levi
subgroup. Write $\frak{l}$ for the Levi algebra associated to $L$.\ We
denote by $P_{L}^{+}$ the set of dominant weights corresponding to $L.$ For
any partition $\lambda \in \mathcal{P}_{n}$ and $\gamma \in P_{L}^{+},$
write $[V^{\frak{g}}(\lambda ):V^{\frak{l}}(\gamma )]$ for the branching
coefficient giving the multiplicity of $V^{\frak{l}}(\gamma )$ (the
irreducible representation of $L$ of highest weight $\gamma $) in the
restriction of $V^{\frak{g}}(\lambda )$ to $L.$

\subsection{Universal characters}

For each Lie algebra $\frak{g=so}_{N}$ or $\frak{sp}_{N}$ and any partition $%
\nu \in \mathcal{P}_{N}$, we denote by $V^{\frak{gl}_{N}}(\nu )\downarrow _{%
\frak{g}}^{\frak{gl}_{N}}$ the restriction of $V^{\frak{gl}_{N}}(\nu )$ to $%
\frak{g}$. Set 
\begin{equation}
V^{\frak{gl}_{N}}(\nu )\downarrow _{\frak{so}_{N}}^{\frak{gl}%
_{N}}=\bigoplus_{\lambda \in \mathcal{P}_{n}}V^{\frak{so}_{N}}(\lambda
)^{\oplus b_{\nu ,\lambda }^{\frak{so}_{N}}}\quad \text{and}\quad V^{\frak{gl%
}_{2n}}(\nu )\downarrow _{\frak{sp}_{2n}}^{\frak{gl}_{2n}}=\bigoplus_{%
\lambda \in \mathcal{P}_{n}}V^{\frak{sp}_{2n}}(\lambda )^{\oplus b_{\nu
,\lambda }^{\frak{sp}_{2n}}}.  \notag
\end{equation}
This makes in particular appear the branching coefficients $b_{\nu ,\lambda
}^{\frak{so}_{N}}$ and $b_{\nu ,\lambda }^{\frak{sp}_{2n}}$.\ The
restriction map $r^{\frak{g}}$ is defined by setting 
\begin{equation*}
r^{\frak{g}}:\left\{ 
\begin{array}{c}
\mathbb{Z[}x_{1},...,x_{N}]^{\mathrm{sym}}\rightarrow \mathcal{R}^{\frak{g}}
\\ 
s_{\nu }^{\frak{gl}_{N}}\longmapsto \mathrm{char}(V^{\frak{gl}_{N}}(\nu
)\downarrow _{\frak{g}}^{\frak{gl}_{N}})
\end{array}
\right. .
\end{equation*}
We have then 
\begin{equation*}
r^{\frak{g}}(s_{\nu }^{\frak{gl}_{N}})=\left\{ 
\begin{array}{c}
s_{\nu }^{\frak{gl}_{N}}(x_{1},...,x_{n},x_{n}^{-1},...,x_{1}^{-1})\text{
when }N=2n \\ 
s_{\nu }^{\frak{gl}_{N}}(x_{1},...,x_{n},1,x_{n}^{-1},...,x_{1}^{-1})\text{
when }N=2n+1
\end{array}
\right. .
\end{equation*}
Let $\mathcal{P}_{n}^{(2)}$ and $\mathcal{P}_{n}^{(1,1)}$ be the subsets of $%
\mathcal{P}_{n}$ containing the partitions with even length rows and the
partitions with even length columns, respectively. When $\nu \in \mathcal{P}%
_{n}$ we have the following formulas for the branching coefficients $b_{\nu
,\lambda }^{\frak{so}_{N}}$ and $b_{\nu ,\lambda }^{\frak{sp}_{2n}}$ :

\begin{proposition}
\label{prop_lit}(see \cite{Li} appendix p\ 295)\label{prop_in_A2n}

\noindent Consider $\nu \in \mathcal{P}_{n}.$ Then:

\begin{enumerate}
\item  $b_{\nu ,\lambda }^{\frak{so}_{2n+1}}=b_{\nu ,\lambda }^{\frak{so}%
_{2n}}=\sum_{\gamma \in \mathcal{P}_{n}^{(2)}}c_{\lambda ,\gamma }^{\nu }$

\item  $b_{\nu ,\lambda }^{\frak{sp}_{2n}}=\sum_{\gamma \in \mathcal{P}%
_{n}^{(1,1)}}c_{\lambda ,\gamma }^{\nu }$
\end{enumerate}

\noindent where $c_{\gamma ,\lambda }^{\nu }$ is the $n$-independent
multiplicity of $s_{\nu }^{\frak{gl}_{n}}$ in the Schur functions product $%
s_{\lambda }^{\frak{gl}_{n}}s_{\gamma }^{\frak{gl}_{n}}$.
\end{proposition}

\noindent\textbf{Remarks: }

\noindent$\mathrm{(i):}$ Note that the equality $b_{\nu,\lambda}^{\frak{so}%
_{2n+1}}=b_{\nu,\lambda}^{\frak{so}_{2n}}$ becomes false in general when $%
\nu\notin\mathcal{P}_{n}$.

\noindent $\mathrm{(ii):}$ By the above proposition we have for any $\nu \in 
\mathcal{P}_{m}$ with $m\leq n$%
\begin{equation}
r^{\frak{sp}_{2n}}(s_{\nu }^{\frak{gl}_{2n}})=\sum_{\lambda \in \mathcal{P}%
_{m}}\sum_{\gamma \in \mathcal{P}_{m}^{(1,1)}}c_{\lambda ,\gamma }^{\nu
}s_{\lambda }^{\frak{sp}_{2n}}\text{ and }r^{\frak{so}_{N}}(s_{\nu }^{\frak{%
gl}_{N}})=\sum_{\lambda \in \mathcal{P}_{m}}\sum_{\gamma \in \mathcal{P}%
_{m}^{(2)}}c_{\lambda ,\gamma }^{\nu }s_{\lambda }^{\frak{so}_{N}}.
\label{dec1}
\end{equation}
By Proposition 1.5.3 in \cite{K}, one has also for any $\lambda \in \mathcal{%
P}_{m}$%
\begin{eqnarray}
s_{\lambda }^{\frak{sp}_{2n}} &=&\sum_{\nu \in \mathcal{P}_{m},\left| \nu
\right| \leq \left| \lambda \right| }(-1)^{\frac{\left| \nu \right| -\left|
\lambda \right| }{2}}\sum_{\alpha =(\alpha _{1}>\cdot \cdot \cdot >\alpha
_{s}>0)}c_{\nu ,\Gamma (\alpha )}^{\lambda }\ r^{\frak{sp}_{2n}}(s_{\nu }^{%
\frak{gl}_{2n}})  \label{dec2} \\
s_{\lambda }^{\frak{so}_{N}} &=&\sum_{\nu \in \mathcal{P}_{m},\left| \nu
\right| \leq \left| \lambda \right| }(-1)^{\frac{\left| \nu \right| -\left|
\lambda \right| }{2}}\sum_{\alpha =(\alpha _{1}>\cdot \cdot \cdot >\alpha
_{s}>0)}c_{\nu ,\Gamma ^{\prime }(\alpha )}^{\lambda }\ r^{\frak{so}%
_{N}}(s_{\nu }^{\frak{gl}_{N}})  \notag
\end{eqnarray}
where $\Gamma (\alpha )=(\alpha _{1}-1,...,\alpha _{s}-1\mid \alpha
_{1},...,\alpha _{s})$ in the Frobenius notation for the partitions. Observe
that the coefficients appearing in the decompositions (\ref{dec1}) and (\ref
{dec2}) do not depend on the rank $n$ considered. Moreover they coincide for
the orthogonal types $B_{n}$ and $D_{n}.$

\bigskip

\noindent As suggested by the above decompositions, the manipulation of the
Weyl characters is simplified by working with infinitely many variables.\ In 
\cite{K}, Koike and Terada have introduced a universal character ring for
the classical Lie groups.\ This ring can be regarded as the ring $\Lambda =%
\mathbb{Z[}x_{1},...,x_{n},...]^{\mathrm{sym}}$ of symmetric functions in
countably many variables.\ It is equipped with three natural $\mathbb{Z}$%
-bases indexed by partitions, namely 
\begin{equation}
\mathcal{B}^{\frak{gl}}\mathcal{=}\{\mathtt{s}_{\lambda }^{\frak{gl}}\mid
\lambda \in \mathcal{P}\},\text{ }\mathcal{B}^{\frak{sp}}\mathcal{=}\{%
\mathtt{s}_{\lambda }^{\frak{sp}}\mid \lambda \in \mathcal{P}\}\text{ and }%
\mathcal{B}^{\frak{so}}\mathcal{=}\{\mathtt{s}_{\lambda }^{\frak{so}}\mid
\lambda \in \mathcal{P}\}.  \label{bases}
\end{equation}
We have then 
\begin{eqnarray}
\mathtt{s}_{\nu }^{\frak{gl}} &=&\sum_{\lambda \in \mathcal{P}}\sum_{\gamma
\in \mathcal{P}^{(2)}}c_{\lambda ,\gamma }^{\nu }\mathtt{s}_{\lambda }^{%
\frak{so}}\text{ and }\mathtt{s}_{\nu }^{\frak{gl}}=\sum_{\lambda \in 
\mathcal{P}}\sum_{\gamma \in \mathcal{P}^{(1,1)}}c_{\lambda ,\gamma }^{\nu }%
\mathtt{s}_{\lambda }^{\frak{sp}}  \label{deuSgl} \\
\mathtt{s}_{\lambda }^{\frak{sp}} &=&\sum_{\nu \in \mathcal{P},\left| \nu
\right| \leq \left| \lambda \right| }(-1)^{\frac{\left| \nu \right| -\left|
\lambda \right| }{2}}\sum_{\alpha =(\alpha _{1}>\cdot \cdot \cdot >\alpha
_{s}>0)}c_{\nu ,\Gamma (\alpha )}^{\lambda }\ \mathtt{s}_{\nu }^{\frak{gl}}
\label{deuSgl2} \\
\mathtt{s}_{\lambda }^{\frak{so}} &=&\sum_{\nu \in \mathcal{P},\left| \nu
\right| \leq \left| \lambda \right| }(-1)^{\frac{\left| \nu \right| -\left|
\lambda \right| }{2}}\sum_{\alpha =(\alpha _{1}>\cdot \cdot \cdot >\alpha
_{s}>0)}c_{\nu ,\Gamma ^{\prime }(\alpha )}^{\lambda }\ \mathtt{s}_{\nu }^{%
\frak{gl}}
\end{eqnarray}
In the sequel we will write for short 
\begin{equation}
b_{\nu ,\lambda }^{\frak{so}}=\sum_{\gamma \in \mathcal{P}^{(2)}}c_{\lambda
,\gamma }^{\nu },\quad b_{\nu ,\lambda }^{\frak{sp}}=\sum_{\gamma \in 
\mathcal{P}^{(1,1)}}c_{\lambda ,\gamma }^{\nu },\quad r_{\lambda ,\nu }^{%
\frak{so}}=\sum_{\alpha }c_{\nu ,\Gamma ^{\prime }(\alpha )}^{\lambda }\text{
and }\quad r_{\lambda ,\nu }^{\frak{sp}}=\sum_{\alpha }c_{\nu ,\Gamma
(\alpha )}^{\lambda }  \label{def_br}
\end{equation}
We denote by $\omega $ the linear involution defined on $\Lambda $ by $%
\omega (\mathtt{s}_{\lambda }^{\frak{gl}})=\mathtt{s}_{\lambda ^{\prime }}^{%
\frak{gl}}.$ Then we have by Theorem 2.3.2 of \cite{K} 
\begin{equation}
\omega (\mathtt{s}_{\lambda }^{\frak{so}})=\mathtt{s}_{\lambda ^{\prime }}^{%
\frak{sp}}.  \label{def_fi}
\end{equation}
Write $\pi _{n}:\mathbb{Z[}x_{1},...,x_{n},...]^{\mathrm{sym}}\rightarrow 
\mathbb{Z[}x_{1},...,x_{n}]^{\mathrm{sym}}$ for the ring homomorphism
obtained by specializing each variable $x_{i},i>n$ at $0.$ Then $\pi _{n}(%
\mathtt{s}_{\lambda }^{\frak{gl}})=s_{\lambda }^{\frak{gl}_{n}}.$ Let $\pi ^{%
\frak{sp}_{2n}}$ and $\pi ^{\frak{so}_{N}}$ be the specialization
homomorphisms defined by setting $\pi ^{\frak{sp}_{2n}}=r^{\frak{sp}%
_{2n}}\circ \pi _{2n}$ and $\pi ^{\frak{so}_{N}}=r^{\frak{so}_{N}}\circ \pi
_{N}$.\ For any partition $\lambda \in \mathcal{P}_{n}$ one has $s_{\lambda
}^{\frak{sp}_{2n}}=\pi ^{\frak{sp}_{2n}}(\mathtt{s}_{\lambda }^{\frak{sp}})$
and $s_{\lambda }^{\frak{so}_{N}}=\pi ^{\frak{so}_{N}}(\mathtt{s}_{\lambda
}^{\frak{s0}}).$ We shall also need the following proposition (see \cite{Ki2}
and \cite{K}).

\begin{proposition}
\label{prop_inde}Consider a Lie algebra $\frak{g}$ of type $X_{n}\in
\{B_{n},C_{n},D_{n}\}$. Let $\lambda \in \mathcal{P}_{r}$ and $\mu \in 
\mathcal{P}_{s}$.\ Suppose $n\geq r+s$ and set 
\begin{equation*}
V^{\frak{g}}(\lambda )\otimes V^{\frak{g}}(\mu )=\bigoplus_{\nu \in \mathcal{%
P}_{n}}V^{\frak{g}}(\nu )^{\oplus d_{\lambda ,\mu }^{\nu }}.
\end{equation*}
Then the coefficients $d_{\lambda ,\mu }^{\nu }$ neither depend on the rank $%
n$ of $\frak{g}$ nor on its type $B,C$ or $D$. More we have 
\begin{equation*}
d_{\lambda ,\mu }^{\nu }=\sum_{\xi ,\sigma ,\tau }c_{\xi ,\sigma }^{\lambda
}c_{\xi ,\tau }^{\mu }c_{\sigma ,\tau }^{\nu }.
\end{equation*}
\end{proposition}

\noindent \textbf{Remarks: }

\noindent $\mathrm{(i):}$ The previous proposition implies the
decompositions $\mathtt{s}_{\lambda }^{\frak{sp}}\times \mathtt{s}_{\mu }^{%
\frak{sp}}=\sum_{\nu \in \mathcal{P}}d_{\lambda ,\mu }^{\nu }\mathtt{s}_{\nu
}^{\frak{sp}}$ and $\mathtt{s}_{\lambda }^{\frak{so}}\times \mathtt{s}_{\mu
}^{\frak{so}}=\sum_{\nu \in \mathcal{P}}d_{\lambda ,\mu }^{\nu }\mathtt{s}%
_{\nu }^{\frak{so}}$ for any $\lambda ,\mu \in \mathcal{P}$, in the ring $%
\Lambda $.

\noindent $\mathrm{(ii):}$ The analogous result for $\frak{g=gl}_{n}$ is
well-known: the outer multiplicities $c_{\lambda ,\mu }^{\nu }$ appearing in
the decomposition of $V^{\frak{gl}_{n}}(\lambda )\otimes V^{\frak{gl}%
_{n}}(\mu )$ do not depend on $n$ provided $n\geq r+s.$

\section{Plethysms and stabilized plethysms\label{sec-stabplethys}}

\subsection{Plethysms on the Weyl characters}

Consider $f\in \Lambda $ and $s_{\lambda }^{\frak{g}}$ the Weyl character
for $\frak{g}$ associated to $\lambda \in \mathcal{P}_{n}$.\ Set $s_{\lambda
}^{\frak{g}}=\sum_{\beta \in \mathbb{Z}^{n}}a_{\beta }x^{\beta }.\;$As in
the case of ordinary plethysms on symmetric functions (see \cite{mac} p
135), one defines the set of variables $y_{i}$ such that 
\begin{equation*}
\prod_{i}(1+ty_{i})=\prod_{\beta }(1+tx^{\beta })^{a_{\beta }}.
\end{equation*}
Then the plethysm of the Weyl character $s_{\lambda }^{\frak{g}}$ by the
symmetric function $f$ is defined by $f\circ s_{\lambda }^{\frak{g}%
}=f(y_{1,}y_{2},...).$ In the sequel, we will focus on the power sum
plethysms $\psi _{\ell }$ where $\ell $ is a positive integer.\ They are
defined from the identity $\psi _{\ell }(s_{\lambda }^{\frak{g}})=p_{\ell
}\circ s_{\lambda }^{\frak{g}}=s_{\lambda }^{\frak{g}}(x_{1}^{\ell
},...,x_{n}^{\ell }).\;$In particular, the map $\psi _{\ell }$ is linear on $%
\mathcal{R}^{\frak{g}}.$ The characters of the symmetric and antisymmetric
parts of $V^{\frak{g}}(\lambda )^{\otimes 2}$ can be expressed as plethysms
by the complete and elementary symmetric functions $h_{2}$ and $e_{2}.$ More
precisely we have 
\begin{equation*}
h_{2}\circ s_{\lambda }^{\frak{g}}=\mathrm{char}(S^{2}(V^{\frak{g}}(\lambda
))\text{ and }e_{2}\circ s_{\lambda }^{\frak{g}}=\mathrm{char}(\Lambda
^{2}(V^{\frak{g}}(\lambda )).
\end{equation*}
From the identities $h_{2}^{2}=\frac{1}{2}(e_{1}^{2}+p_{2})$ and $e^{2}=%
\frac{1}{2}(e_{1}^{2}-p_{2})$, we derive the relations 
\begin{equation}
h_{2}\circ s_{\lambda }^{\frak{g}}=\frac{1}{2}((s_{\lambda }^{\frak{g}%
})^{2}+\psi _{2}(s_{\lambda }^{\frak{g}}))\text{ and }e_{2}\circ s_{\lambda
}^{\frak{g}}=\frac{1}{2}((s_{\lambda }^{\frak{g}})^{2}-\psi _{2}(s_{\lambda
}^{\frak{g}})).  \label{relfund}
\end{equation}

\subsection{Stabilized plethysms on the Schur functions\label%
{subsec_ple_shur}}

Given $(\mu ^{(0)},...,\mu ^{(\ell -1})$ a $\ell $-tuple of partitions, we
write $c_{(\mu ^{(0)},...,\mu ^{(\ell -1})}^{\lambda }$ for the $n$%
-independent coefficient of $s_{\lambda }^{\frak{gl}_{n}}$ in the product $%
s_{\mu ^{(0)}}^{\frak{gl}_{n}}\cdots s_{\mu ^{(\ell -1)}}^{\frak{gl}_{n}}$.\
For any partition $\lambda \in \mathcal{P}_{n},$ the plethysm $\psi _{\ell
}(s_{\lambda }^{\frak{gl}_{n}})$ decomposes on the basis of Schur functions
on the form 
\begin{equation}
\psi _{\ell }(s_{\lambda }^{\frak{gl}_{n}})=\sum_{\left| \mu \right| =\ell
\left| \lambda \right| }\varepsilon (\mu )c_{(\mu ^{(0)},...,\mu ^{(\ell
-1)})}^{\lambda }s_{\mu }^{\frak{gl}_{n}}.  \label{decA}
\end{equation}
Here $\varepsilon (\mu )\in \{-1,0,1\}$ and $\mu /\ell =(\mu ^{(0)},...,\mu
^{(\ell -1})$ are respectively the $\ell $-sign and the $\ell $-quotient of
the partition $\mu $.\ We now briefly recall the algorithm which permits to
obtain the sign $\varepsilon (\mu )$ and the $\ell $-tuple of partitions $%
\mu /\ell .\;$Our description slightly differs from that which can be
usually found in the literature (see \cite{mac} Example 8 p 12).\ This is
because we have made our notation consistent with Section \ref{sec-algo}.

\noindent Set $\rho _{n}=(1,2,...,n)$ and $I_{n}=\{1,2,...,n\}.\;$For any $%
k\in \{0,...,\ell -1\}$ consider the sequences 
\begin{equation*}
I^{(k)}=(i\in I_{n}\mid \mu _{i}+i\equiv k\func{mod}\ell )\text{ and }%
J^{(k)}=(i\in I_{n}\mid i\equiv k\func{mod}\ell )
\end{equation*}
in which the entries occur in the increasing order.\ Set $r_{k}=\mathrm{card}%
(I^{(k)})$ and write $I^{(k)}=(i_{1}^{(k)},...,i_{r_{k}}^{(k)})$.

\begin{enumerate}
\item  If there exists $k\in \{0,...,\ell -1\}$ such that $\mathrm{card}%
(I^{(k)})\neq \mathrm{card}(J^{(k)})$ then $\varepsilon (\mu )=0.$

\item  Otherwise let $\sigma _{0}\in S_{n}$ be the permutation mapping $%
I^{(k)}$ to $J^{(k)}$ for any $k=0,...,\ell -1$.\ Then we have $\varepsilon
(\mu )=\varepsilon (\sigma _{0})$ and $\mu /\ell =(\mu ^{(0)},...,\mu
^{(\ell -1})$ where for any $k\in \{0,...,\ell -1\}$%
\begin{equation}
\mu ^{(k)}=\left( \frac{\mu _{i}+i+\ell -k}{\ell }\mid i\in I^{(k)}\right)
-(1,2,...,r_{k})\in \mathbb{Z}^{r_{k}}  \label{def_mukA}
\end{equation}
(see example below).
\end{enumerate}

\noindent \textbf{Remark: }Set $n=q\ell +r$ where $q$ and $r$ are
respectively the quotient and the rest of the division of $n$ by $\ell $.
Then we have $\mathrm{card}(J^{(k)})=q+1$ for any $k\in \{1,\ldots ,r\}$ and 
$\mathrm{card}(J^{(k)})=q$ for any $k\in \{0,r+1,\ldots ,n\}$. Hence in (\ref
{def_mukA}), we have $r_{k}=q+1$ for any $k\in \{1,\ldots ,r\}$ and $r_{k}=q$
for any $k\in \{0,r+1,\ldots ,n\}$ (see also Remark just before Section 4.2).

\begin{example}
\ \ \label{ex_A}

\noindent Consider $\mu =(1,2,3,4,4,4,6,6)$ and take $\ell =3.$ We have $\mu
+\rho _{8}=(2,4,6,8,9,10,13,14).$ Thus $%
I^{(0)}=(3,5),I^{(1)}=(2,6,7),I^{(2)}=(1,4,8)$ and $%
J^{(0)}=(3,6),J^{(1)}=(1,4,7),J^{(2)}=(2,5,8).$ Then $\mu ^{(0)}=(1,1),\mu
^{(1)}=(1,2,2)$ and $\mu ^{(2)}=(0,1,2).$ Moreover 
\begin{equation*}
\sigma _{0}=\left( 
\begin{array}{cccccccc}
1 & 2 & 3 & 4 & 5 & 6 & 7 & 8 \\ 
2 & 1 & 3 & 5 & 6 & 4 & 7 & 8
\end{array}
\right) .
\end{equation*}
Hence $\varepsilon (\mu )=-1$.
\end{example}

\begin{proposition}
\label{prop_stabA}Consider $\mu \in \mathcal{P}_{n}$ such that $\varepsilon
(\mu )\neq 0$ and set $\mu /\ell =(\mu ^{(0)},...,\mu ^{(\ell -1}).\;$Let $%
\nu \in \mathcal{P}_{n+1}$ be the partition obtained by adding in $\mu $ a
part $0$.\ Then $\varepsilon (\nu )=\varepsilon (\mu )$ and we have $\nu
/\ell =(\nu ^{(0)},...,\nu ^{(\ell -1})$ where $\nu ^{(0)}=\mu ^{(\ell -1)},$
$\nu ^{(k)}=\mu ^{(k-1)}$ for any $k\in \{2,\ldots ,\ell -1\}$ and $\nu
^{(1)}=(0,\mu ^{(0)})$ is obtained by adding a part $0$ in $\mu ^{(0)}$.
\end{proposition}

\begin{proof}
Let us slightly abuse the notation and write $I^{(k)}(\mu ),$ $J^{(k)}(\mu ),
$ $I^{(k)}(\nu ),$ $J^{(k)}(\nu ),$ $k=0,\ldots ,\ell -1$ for the sequences
defined from $\mu $ and $\nu $ by applying the previous procedure. Then, we
have 
\begin{equation}
\left\{ 
\begin{tabular}{l}
$I^{(1)}(\nu )=\{1\}\cup (I^{(0)}(\mu )+1),I^{(0)}(\nu )=(I^{(\ell -1)}(\mu
)+1)$ \\ 
$I^{(k)}(\nu )=(I^{(k-1)}(\mu )+1)\text{ for }k=2,\ldots \ell -1.$%
\end{tabular}
\right.   \label{relI}
\end{equation}
Here by $(I^{(k-1)}(\mu )+1),$ we mean the sequence obtained by adding $1$
to the entries of $I^{(k-1)}(\mu )$.\ Set $n=q\ell +r$ as in the previous
remark. We will assume that $r\neq \ell -1$ so that $q$ and $r+1$ are
respectively the quotient and the rest of the division of $n+1$ by $\ell $.\
The case $r=\ell -1$ is similar.\ We have then $\mathrm{card}(J^{(k)}(\mu ))=%
\mathrm{card}(I^{(k)}(\mu ))=q+1$ for $k\in \{1,\ldots ,r\}$ and $\mathrm{%
card}(J^{(k)}(\mu ))=\mathrm{card}(I^{(k)}(\mu ))=q$ for $k\in
\{0,r+1,\ldots ,n\}$. Now observe that $J^{(k)}(\nu )=J^{(k)}(\mu )$ for $%
k\neq r+1$ and $J^{(r+1)}(\nu )=J^{(r+1)}(\mu )\cup \{n+1\}$.\ This implies
that $\mathrm{card}(J^{(k)}(\nu ))=\mathrm{card}(I^{(k)}(\nu ))=q+1$ for $%
k\in \{1,\ldots ,r+1\}$ and $\mathrm{card}(J^{(k)}(\mu ))=\mathrm{card}%
(I^{(k)}(\mu ))=q$ for $k\in \{0,r+2,\ldots ,n\}$. Thus $\varepsilon (\nu
)=\varepsilon (\mu )$. We then easily deduce $\nu ^{(0)},...,\nu ^{(\ell -1}$
from (\ref{def_mukA}) and (\ref{relI}).
\end{proof}

\bigskip

\noindent \textbf{Remarks: }

\noindent $\mathrm{(i):}$ The the decomposition (\ref{decA}) do not depend
on the rank $n$ considered provided $n>\ell l(\lambda )$.\ Indeed, by
Proposition \ref{prop_stabA}, $\varepsilon (\mu )$ and the non-zero parts of
the partitions $\mu ^{(k)}$ of the above algorithm are not modified when
empty parts are added to $\mu $.

\noindent $\mathrm{(ii):}$ When $n>\ell l(\lambda )$, we write for short $%
a_{\lambda ,\mu }^{\ell ,\frak{gl}}=\varepsilon (\mu )c_{(\mu ^{(0)},...,\mu
^{(\ell -1})}^{\lambda }.$ Then $a_{\lambda ,\mu }^{\ell ,\frak{gl}}\neq 0$
only if $\left| \mu \right| =\ell \left| \lambda \right| .$

\begin{proposition}
\label{stab_pleth_schur}Consider $f\in \Lambda $ with degree $d$ and $%
\lambda \in \mathcal{P}_{n}.\;$Then the coefficients of the expansion of $%
f\circ s_{\lambda }^{\frak{gl}_{n}}$ on the basis of Schur functions do not
depend on $n$ provided $n>dl(\lambda )$.
\end{proposition}

\begin{proof}
By (\ref{decA}) and the previous remark, the proposition is true for the
power sum plethysms $p_{\ell }\circ s_{\lambda }^{\frak{gl}_{n}}.\;$The map $%
g\mapsto g\circ s_{\lambda }^{\frak{gl}_{n}}$ is a ring homomorphism of $%
\Lambda _{n}.\;$The subspace $\Lambda _{n}^{d}$ of polynomials in $\Lambda
_{n}$ with degree $d$ is generated by the Newton polynomials $p_{\beta
}=p_{\beta _{1}}\cdot \cdot \cdot p_{\beta _{k}}$, such that $\beta _{i}\in 
\mathbb{N}$ and $\beta _{1}+\cdot \cdot \cdot +\beta _{k}=d.\;$So it
suffices to prove the proposition for $f=p_{\beta }.$ We have $p_{\beta
}\circ s_{\lambda }^{\frak{gl}_{n}}=p_{\beta _{1}}\circ s_{\lambda }^{\frak{%
gl}_{n}}\times \cdot \cdot \cdot \times p_{\beta _{k}}\circ s_{\lambda }^{%
\frak{gl}_{n}}.\;$Suppose $n>dl(\lambda ).\;$For any $i=1,...,k,$ we have $%
n\geq \beta _{i}l(\lambda ).\;$Thus we can write $p_{\beta _{i}}\circ
s_{\lambda }^{\frak{gl}_{n}}=\sum_{\mu ^{(i)}\in \mathcal{P}_{n}}a_{\lambda
,\mu ^{(i)}}^{\beta _{i},\frak{gl}}\,s_{\mu ^{(i)}}^{\frak{gl}_{n}}.\;$%
Moreover $a_{\lambda ,\mu ^{(i)}}^{\beta _{i},\frak{gl}}\neq 0$ only if $%
\left| \mu ^{(i)}\right| =\beta _{i}\left| \lambda \right| .\;$By Remark $%
\mathrm{(ii)}$ following Proposition \ref{prop_inde}, we obtain that the
coefficients of the decomposition of $p_{\beta }\circ s_{\lambda }^{\frak{gl}%
_{n}}$ on the basis of Schur functions do not depend on $n$ when $n\geq
\left| \beta \right| l(\lambda ).$
\end{proof}

\subsection{Stabilized plethysms on the Weyl characters\label{subsec-stab}}

\begin{lemma}
\label{lem_stab-pl}Consider $\lambda \in \mathcal{P}_{m},$ $\ell $ a
positive integer and $\frak{g}$ an orthogonal or symplectic Lie algebra with
rank $n\geq m.$

\begin{itemize}
\item  The coefficients of the expansion of the plethysm $p_{\ell }\circ
s_{\lambda }^{\frak{g}}$ on the basis of Weyl characters do not depend on $n$
provided $n>\ell l(\lambda ).$

\item  In this case, these coefficients coincide for $\frak{g=so}_{2n+1}$
and $\frak{g=so}_{2n}.$

\item  For any $n>\ell l(\lambda )$, set 
\begin{equation*}
p_{\ell }\circ s_{\lambda }^{\frak{so}_{N}}=\sum_{\mu \in \mathcal{P}%
_{n}}a_{\lambda ,\mu }^{\ell ,\frak{so}}s_{\mu }^{\frak{so}_{N}}\text{ and }%
p_{\ell }\circ s_{\lambda }^{\frak{sp}_{2n}}=\sum_{\mu \in \mathcal{P}%
_{n}}a_{\lambda ,\mu }^{\ell ,\frak{sp}}s_{\mu }^{\frak{sp}_{2n}}.
\end{equation*}
We have 
\begin{eqnarray*}
a_{\lambda ,\mu }^{\ell ,\frak{so}} &=&\sum_{\nu \in \mathcal{P}_{m},\left|
\nu \right| \leq \left| \lambda \right| }(-1)^{\frac{\left| \lambda \right|
-\left| \nu \right| }{2}}\sum_{\delta \in \mathcal{P},\left| \delta \right|
=\ell \left| \nu \right| }r_{\lambda ,\nu }^{\frak{so}}\ a_{\nu ,\delta
}^{\ell ,\frak{gl}}\ b_{\delta ,\mu }^{\frak{so}}, \\
a_{\lambda ,\mu }^{\ell ,\frak{sp}} &=&\sum_{\nu \in \mathcal{P}_{m},\left|
\nu \right| \leq \left| \lambda \right| }(-1)^{\frac{\left| \lambda \right|
-\left| \nu \right| }{2}}\sum_{\delta \in \mathcal{P},\left| \delta \right|
=\ell \left| \nu \right| }r_{\lambda ,\nu }^{\frak{sp}}\ a_{\nu ,\delta
}^{\ell ,\frak{gl}}\ b_{\delta ,\mu }^{\frak{sp}}.
\end{eqnarray*}
\end{itemize}
\end{lemma}

\begin{proof}
We have $n>\ell l(\lambda ).$ Hence, the decomposition $s_{\lambda }^{\frak{%
so}_{N}}=\sum_{\nu \in \mathcal{P}_{m},\left| \nu \right| \leq \left|
\lambda \right| }(-1)^{\frac{\left| \lambda \right| -\left| \nu \right| }{2}%
}r_{\lambda ,\nu }^{\frak{so}}r^{\frak{so}_{N}}(s_{\nu }^{\frak{gl}_{N}})$
holds.\ Since $\psi _{\ell }$ and $r^{\frak{so}_{N}}$ commute, this gives 
\begin{multline*}
p_{\ell }\circ s_{\lambda }^{\frak{so}_{N}}=\sum_{\nu \in \mathcal{P}%
_{m},\left| \nu \right| \leq \left| \lambda \right| }(-1)^{\frac{\left|
\lambda \right| -\left| \nu \right| }{2}}\sum_{\delta \in \mathcal{P},\left|
\delta \right| =\ell \left| \nu \right| }r_{\lambda ,\nu }^{\frak{so}}\
a_{\nu ,\delta }^{\ell ,\frak{gl}}\,r^{\frak{so}_{N}}(s_{\delta }^{\frak{gl}%
_{N}})= \\
\sum_{\mu \in \mathcal{P}_{n}}\sum_{\nu \in \mathcal{P}_{m},\left| \nu
\right| \leq \left| \lambda \right| }(-1)^{\frac{\left| \lambda \right|
-\left| \nu \right| }{2}}\sum_{\delta \in \mathcal{P},\left| \delta \right|
=\ell \left| \nu \right| }r_{\lambda ,\nu }^{\frak{so}}\ a_{\nu ,\delta
}^{\ell ,\frak{gl}}\ b_{\delta ,\mu }^{\frak{so}}\,s_{\mu }^{\frak{so}_{N}}.
\end{multline*}
This yields the desired expression for the coefficients $a_{\lambda ,\mu
}^{\ell ,\frak{so}}.\;$In particular they do not depend on $n$ and coincide
for $\frak{g=so}_{2n+1}$ and $\frak{g=so}_{2n}.$ The proof is similar for $%
\frak{g=sp}_{2n}.$
\end{proof}

\begin{proposition}
\label{prop-stabpleth}Consider $f\in \Lambda $ with degree $d$ and $\lambda
\in \mathcal{P}_{n}.\;$Then the coefficients of the expansion of $f\circ
s_{\lambda }^{\frak{g}}$ on the basis of Schur functions do not depend on $n$
provided $n\geq d\left| \lambda \right| $. In this case, these coefficients
coincide for $\frak{g=so}_{2n+1}$ and $\frak{g=so}_{2n}.$
\end{proposition}

\begin{proof}
The proposition follows from Lemma \ref{lem_stab-pl} by similar arguments to
those of Proposition \ref{stab_pleth_schur}.
\end{proof}

\bigskip

\noindent According to the previous Lemma, we have the decompositions 
\begin{equation*}
p_{\ell }\circ \mathtt{s}_{\lambda }^{\frak{gl}}=\sum_{\mu }a_{\lambda ,\mu
}^{\ell ,\frak{gl}}\mathtt{s}_{\mu }^{\frak{gl}},\quad p_{\ell }\circ 
\mathtt{s}_{\lambda }^{\frak{sp}}=\sum_{\mu }a_{\lambda ,\mu }^{\ell ,\frak{%
sp}}\mathtt{s}_{\mu }^{\frak{sp}}\text{ and }\quad p_{\ell }\circ \mathtt{s}%
_{\lambda }^{\frak{so}}=\sum_{\mu }a_{\lambda ,\mu }^{\ell ,\frak{so}}%
\mathtt{s}_{\mu }^{\frak{so}}.
\end{equation*}
We shall need in Section \ref{subsec_main-TH} the following Lemma :

\begin{lemma}
\label{lem_orthsympl_dual}Consider $f\in \Lambda $ and $\lambda \in \mathcal{%
P}.$ Then

\begin{itemize}
\item  $\omega (f\circ \mathtt{s}_{\lambda }^{\frak{g}})=f\circ \omega (%
\mathtt{s}_{\lambda }^{\frak{g}})$ if $\left| \lambda \right| $ is even,

\item  $\omega (f\circ \mathtt{s}_{\lambda }^{\frak{g}})=\omega (f)\circ
\omega (\mathtt{s}_{\lambda }^{\frak{g}})$ if $\left| \lambda \right| $ is
odd.
\end{itemize}
\end{lemma}

\begin{proof}
From Example 1 page 136 of \cite{mac} we have for any positive integer $\ell 
$, $\omega (p_{\ell }\circ g)=p_{\ell }\circ \omega (g)$ if $g$ is
homogeneous of even degree and $\omega (p_{\ell }\circ g)=\omega (p_{\ell
})\circ \omega (g)$ if $g$ is homogeneous of odd degree.\ Since $\psi _{\ell
}$ is linear, this shows that $\omega (p_{\ell }\circ \mathtt{s}_{\lambda }^{%
\frak{g}})=p_{\ell }\circ \omega (\mathtt{s}_{\lambda }^{\frak{g}})$ if $%
\left| \lambda \right| $ is even and $\omega (p_{\ell }\circ \mathtt{s}%
_{\lambda }^{\frak{g}})=\omega (p_{\ell })\circ \omega (\mathtt{s}_{\lambda
}^{\frak{g}})$ if $\left| \lambda \right| $ is odd. Indeed, according to (%
\ref{deuSgl2}), $\mathtt{s}_{\lambda }^{\frak{g}}$ is a sum of homogeneous
functions of degrees equal to $\left| \lambda \right| $ modulo $2$.\ The
Lemma then follows since the maps $\omega $ and $f\mapsto f\circ \mathtt{s}%
_{\lambda }^{\frak{g}}$ are ring homomorphisms of $\Lambda .$
\end{proof}

\bigskip

\noindent \textbf{Remarks:}

\noindent $\mathrm{(i):}$ Since $\omega (p_{\ell })=(-1)^{\ell -1}p_{\ell },$
one has by the previous lemma $a_{\lambda ,\mu }^{\ell ,\frak{sp}%
}=a_{\lambda ^{\prime },\mu ^{\prime }}^{\ell ,\frak{so}}$ if $\left|
\lambda \right| $ is even and $a_{\lambda ,\mu }^{\ell ,\frak{sp}%
}=(-1)^{\ell -1}a_{\lambda ^{\prime },\mu ^{\prime }}^{\ell ,\frak{so}}$
otherwise. This can also be verified by using the explicit formulas of Lemma 
\ref{lem_stab-pl}.

\noindent $\mathrm{(ii):}$ The coefficients $a_{\lambda ,\mu }^{\ell ,\frak{%
so}}$ are rather complicated to compute by using formulas of Lemma \ref
{lem_stab-pl}. We are going to see in the following Section that they
coincide with branching coefficients corresponding to restriction to certain
Levi subalgebras.

\section{Power sum plethysms for Weyl characters of type $B_{n}\label%
{sec-algo}$}

\subsection{Statement of the theorem}

In Theorem 3.2.8 of \cite{lec}, we have described an algorithm for computing
the plethysms $p_{\ell }\circ s_{\lambda }^{\frak{so}_{2n+1}}$ for any
positive integer $\ell $ and any rank $n.\;$It notably permits to show that
the decomposition of $p_{\ell }\circ s_{\lambda }^{\frak{so}_{2n+1}}$ on the
basis of Weyl characters makes appear branching coefficients corresponding
to the restriction to a Levi subgroup of $\frak{so}_{2n+1}.$ Surprisingly,
similar algorithms for $\frak{sp}_{2n}$ and $\frak{so}_{2n}$ only exists
when $\ell $ is odd.\ In particular, the coefficients of the decomposition
of $p_{\ell }\circ s_{\lambda }^{\frak{sp}_{2n}}$ and $p_{\ell }\circ
s_{\lambda }^{\frak{so}_{2n}}$ on the basis of Weyl characters are not
branching coefficients in general when $\ell $ is even. As we are going to
see, this is nevertheless the case for the stabilized forms of these
plethysms.

\bigskip

Theorems 3.2.8 and 3.2.10 of \cite{lec} can be reformulated as follows :

\begin{theorem}
\label{Th-decB}For any partition $\lambda \in \mathcal{P}_{n}$ and any
positive integer $\ell $ we have 
\begin{equation}
p_{\ell }\circ s_{\lambda }^{\frak{so}_{2n+1}}=\sum_{\mu \in \mathcal{P}%
_{n}}\varepsilon (\mu )[V^{\frak{so}_{2n+1}}(\lambda ):V^{\frak{g}_{\ell
,\mu }}(\gamma _{\ell ,\mu })]s_{\mu }^{\frak{so}_{2n+1}}  \label{dec6levi}
\end{equation}
where

\begin{itemize}
\item  $\varepsilon (\mu )\in \{-1,0,1\}$,

\item  $\frak{g}_{\ell ,\mu }$\ is the Levi algebra of $G_{\ell ,\mu }$, a
Levi subgroup of $SO_{2n+1},$

\item  $\gamma _{\ell ,\mu }$ is a dominant weight for $G_{\ell ,\mu }$.
\end{itemize}

\noindent Moreover, $\varepsilon (\mu ),$ $G_{\ell ,\mu }$ and $\gamma
_{\ell ,\mu }$ are determined from $\mu $ and $\ell $ by an algorithm which
can be regarded as an analogue in type $B_{n}$ of the computation of the $%
\ell $-quotient $\mu /\ell $.
\end{theorem}

\noindent We now recall the algorithm which permits to determinate $%
\varepsilon (\mu ),$ $G_{\ell ,\mu }$ and $\gamma _{\ell ,\mu }$ in the
above theorem. Set 
\begin{equation*}
J_{n}=\{\overline{n},...,\overline{1},1,...,n\}\text{ and }L_{n}=\{\overline{%
n-1},...,\overline{1},0,1,...,n\}.
\end{equation*}
Let $\eta $ be the bijection from $J_{n}$ to $L_{n}$ defined by $\eta
(x)=x+1 $ if $x<0$ and $\eta (x)=x$ otherwise. For each element $w\in W$
(the Weyl group of $\frak{so}_{2n+1}$), denote by $\widetilde{w}$ the
bijection from $J_{n}$ to $L_{n}$ defined by $\widetilde{w}=\eta \circ w.$
This means that $\widetilde{w}(x)=w(x)$ if $w(x)>0$ and $\widetilde{w}%
(x)=w(x)+1$ if $w(x)<0.\;$In particular $w$ is determined by $\widetilde{w}%
.\;$For any $x\in L_{n},$ set $x^{\ast }=\overline{x}+1.$ The map $%
x\longmapsto x^{\ast }$ is involutive from $L_{n}$ to itself.$\;$Since $w(%
\overline{x})=\overline{w(x)}$, we have also 
\begin{equation}
\widetilde{w}(\overline{x})=\widetilde{w}(x)^{\ast }.  \label{rel}
\end{equation}
Hence, $\widetilde{w}$ is determined by the images of any subset $%
U_{n}\subset J_{n}$ such that $\mathrm{card}(U_{n})=n$ and $x\in U_{n}$
implies $\overline{x}\notin U_{n}.$

\noindent For any $k=1,...,\ell $ set 
\begin{equation}
I^{(k)}=(i\in I_{n}\mid \mu _{i}+i\equiv k\func{mod}\ell )\text{ and }%
J^{(k)}=(x\in L_{n}\mid x\equiv k\func{mod}\ell ).  \label{IJB}
\end{equation}
Note that $(J^{(k)})^{\ast }=J^{(l-k+1)}$.

\bigskip 

\noindent \textbf{Remark: }Set $n=q\ell +r$ where $q$ and $r$ are
respectively the quotient and the rest of the division of $n$ by $\ell $.
Then we have 
\begin{eqnarray}
\mathrm{card}(J^{(k)}) &=&\left\{ 
\begin{array}{l}
2q\text{ for }\mathrm{min}(r+1,\ell -r-1)\leq k\leq \mathrm{max}(r+1,\ell
-r-1) \\ 
2q+1\text{ otherwise}
\end{array}
\right. \text{ when }r\neq \frac{\ell }{2},  \label{carJk} \\
\mathrm{card}(J^{(k)}) &=&2q+1\text{ for any }k\in \{1,\ldots ,\ell \}\text{
when }r=\frac{\ell }{2}.  \notag
\end{eqnarray}

\subsection{The even case $\ell =2p\label{subsec-algoeven}$}

For any $k=1,...,p$, set $s_{k}=\mathrm{card}(I^{(k)}),$ $r_{k}=\mathrm{card}%
(I^{(k)})+\mathrm{card}(I^{(\ell -k+1)})$ and define $X^{(k)}$ as the
increasing reordering of $\overline{I}^{(k)}\cup I^{(\ell -k+1)}.\;$Set 
\begin{equation}
X^{(k)}=(i_{1}^{(k)},...,i_{r_{k}}^{(k)}).  \label{defXKB}
\end{equation}

\begin{enumerate}
\item  If there exists $k\in \{1,...,p\}$ such that $\mathrm{card}%
(X^{(k)})\neq \mathrm{card}(J^{(k)})$ then $\varepsilon (\mu )=0.$

\item  Otherwise we have $\mathrm{card}(J^{(\ell -k+1)})=\mathrm{card}%
(J^{(k)})=r_{k}$ since $(J^{(k)})^{\ast }=J^{(l-k+1)}.\;$Let $w_{0}$ be the
unique element of $W$ mapping $X^{(k)}$ to $J^{(l-k+1)}$ for any $k=1,...,p$%
. Define $\alpha _{k}=\frac{1}{\ell }(\max J^{(k)}-k)$. For any $k=1,...,p$,
consider $\mu ^{(k)}\in \widetilde{\mathcal{P}}_{r_{k}}$ defined by 
\begin{equation}
\mu ^{(k)}=\left( \mathrm{sign}(i)\frac{\mu _{\left| i\right| }+\left|
i\right| +\mathrm{sign}(i)k-\frac{1+\mathrm{sign}(i)}{2}}{\ell }\mid i\in
X^{(k)}\right) -(1,...,r_{k})+(\alpha _{k}+1,....,\alpha _{k}+1).
\label{mukB}
\end{equation}
\bigskip
\end{enumerate}

\noindent \textbf{Remark: }With $q$ and $r$ as in (\ref{carJk}), we can
write 
\begin{eqnarray*}
\alpha _{k} &=&\left\{ 
\begin{tabular}{l}
$q-1\text{ for }\mathrm{min}(r+1,\ell -r-1)\leq k\leq \mathrm{max}(r+1,\ell
-r-1)$ \\ 
$q$ otherwise
\end{tabular}
\right. \text{ when }r\neq \frac{\ell }{2}, \\
\alpha _{k} &=&q\text{ for any }k\in \{1,\ldots ,\ell \}\text{ when }r=\frac{%
\ell }{2}.
\end{eqnarray*}
Note also that in step 2, $r_{k}\in \{2q,2q+1\}$ according to (\ref{carJk}).

\bigskip

We have then with the above notation : 
\begin{equation*}
\varepsilon (\mu )=\varepsilon (w_{0}),\quad G_{\ell ,\mu }=GL_{r_{1}}\times
\cdots \times GL_{r_{p}}\quad \text{and }\quad \gamma _{\ell ,\mu }=(\mu
^{(1)},...,\mu ^{(p)})\in P_{G_{\ell ,\mu }}^{+}.
\end{equation*}

\begin{example}
\label{ex-2}Put $n=6,$ $\ell =2$ (thus $p=1$) and consider $\mu
=(2,5,5,6,7,9).\;$Then $\mu +\rho _{6}=(3,7,8,10,12,15).\;$Hence $%
I^{(2)}=(3,4,5)$ and $I^{(1)}=(1,2,6).\;$Moreover $J^{(1)}=(\overline{5},%
\overline{3},\overline{1},1,3,5)$ and $J^{(2)}=(\overline{4},\overline{2}%
,0,2,4,6).$ Then $\widetilde{w}_{0}$ sends $X^{(1)}=(\overline{6},\overline{2%
},\overline{1},3,4,5)$ on $J^{(2)}.$ This gives 
\begin{equation*}
\widetilde{w}_{0}=\left( 
\begin{array}{cccccccccccc}
\overline{6} & \overline{5} & \overline{4} & \overline{3} & \overline{2} & 
\overline{1} & 1 & 2 & 3 & 4 & 5 & 6 \\ 
\overline{4} & \overline{5} & \overline{3} & \overline{1} & \overline{2} & 0
& 1 & 3 & 2 & 4 & 6 & 5
\end{array}
\right)
\end{equation*}
by using (\ref{rel}).\ Hence 
\begin{equation*}
w_{0}=\left( 
\begin{array}{cccccccccccc}
\overline{6} & \overline{5} & \overline{4} & \overline{3} & \overline{2} & 
\overline{1} & 1 & 2 & 3 & 4 & 5 & 6 \\ 
\overline{5} & \overline{6} & \overline{4} & \overline{2} & \overline{3} & 
\overline{1} & 1 & 3 & 2 & 4 & 6 & 5
\end{array}
\right) .
\end{equation*}
We have $\varepsilon (\mu )=1,$ $\alpha _{1}=2$ and $\gamma _{\ell ,\mu
}=(\mu ^{(1)})$ where 
\begin{equation*}
\mu ^{(1)}=(-7,-3,-1,4,5,6)-(1,2,3,4,5,6)+(3,3,3,3,3,3)=(-5,-2,-1,3,3,3).
\end{equation*}
Observe that $G_{\ell ,\mu }\simeq GL_{6}$.
\end{example}

\subsection{The odd case $\ell =2p+1\label{subsec-algoodd}$}

In addition to the sets $X^{(k)},k=1,...,p$ defined in (\ref{defXKB}), we
have also to consider $I^{(p+1)}$.$\;$Set $r_{p+1}=\mathrm{card}(I^{(p+1)})$
and write $I^{(p+1)}=\{i_{1}^{(p+1)},...,i_{r_{p+1}}^{(p+1)}\}.$ Observe
that $(J^{(p+1)})^{\ast }=J^{(p+1)}$. Let $X^{(p+1)}$ be the increasing
reordering of $\overline{I}^{(p+1)}\cup I^{(p+1)}$.

\begin{enumerate}
\item  If $\mathrm{card}(I^{(p+1)})\neq \frac{1}{2}\mathrm{card}(J^{(p+1)})$
or if there exists $k\in \{1,...,p\}$ such that $\mathrm{card}(X^{(k)})\neq 
\mathrm{card}(J^{(k)})$ then $\varepsilon (\mu )=0.$

\item  Otherwise, we have $\mathrm{card}(J^{(p+1)})=2\mathrm{card}%
(I^{(p+1)})=2r_{p+1}$. Let $w_{0}$ be the unique element of $W$ mapping $%
X^{(k)}$ to $J^{(l-k+1)}$ for any $k=1,...,p$ and $X^{(p+1)}$ to $%
J^{(p+1)}.\;$Define 
\begin{equation*}
\mu ^{(p+1)}=\left( \frac{\mu _{i}+i+p}{\ell }\mid i\in I^{(p+1)}\right)
-(1,...,r_{p+1})\in \mathcal{P}_{r_{p+1}}
\end{equation*}
and for any $k=1,...,p,$ $\mu ^{(k)}$ as in the even case. Set $\mathcal{I=}%
\{I^{(p+1)},X^{(1)},...,X^{(p)}\}.$ We have then with the above notation 
\begin{equation*}
\varepsilon (\mu )=\varepsilon (w_{0}),\quad G_{\ell ,\mu }=GL_{r_{1}}\times
\cdots \times GL_{r_{p}}\times GL_{2r_{p+1}+1}\quad \text{and}\quad \gamma
_{\ell ,\mu }=(\mu ^{(p+1)},\mu ^{(1)},...,\mu ^{(p)})\in P_{G_{\ell ,\mu
}}^{+}.
\end{equation*}
\end{enumerate}

\noindent \textbf{Remark: }With $q$ and $r$ as in (\ref{carJk}), we have $%
r_{p+1}=\frac{1}{2}\mathrm{card}(J^{(p+1)})=q$ when $1$ is satisfied.

\begin{example}
Put $n=6,$ $\ell =3$ (thus $p=1$) and consider $\mu =(1,5,5,6,7,9).$ We have 
$\mu +\rho _{6}=(2,7,8,10,12,15).$ Thus $X^{(1)}=(\overline{4},\overline{2}%
,5,6),I^{(1)}=(2,4),I^{(2)}=(1,3)$ and $J^{(1)}=(\overline{5},\overline{2}%
,1,4),J^{(2)}=(\overline{4},\overline{1},2,5)$ and $J^{(3)}=(\overline{3}%
,0,3,6)$. In particular $\alpha _{1}=1.\;$Then 
\begin{equation*}
\mu ^{(1)}=\left( -\frac{10-1}{3}-1+2,-\frac{7-1}{3}-2+2,\frac{12}{3}-3+2,%
\frac{15}{3}-4+2\right) =(-2,-2,3,3)
\end{equation*}
and $\mu ^{(2)}=(\frac{2+1}{3}-1,\frac{8+1}{3}-2)=(0,1).$ Moreover, one has
by using (\ref{rel}) 
\begin{equation*}
\widetilde{w}_{0}=\left( 
\begin{array}{cccccccccccc}
\overline{6} & \overline{5} & \overline{4} & \overline{3} & \overline{2} & 
\overline{1} & 1 & 2 & 3 & 4 & 5 & 6 \\ 
\overline{5} & \overline{2} & \overline{3} & \overline{4} & 0 & \overline{1}
& 2 & 1 & 5 & 4 & 3 & 6
\end{array}
\right) .
\end{equation*}
Hence 
\begin{equation*}
w_{0}=\left( 
\begin{array}{cccccccccccc}
\overline{6} & \overline{5} & \overline{4} & \overline{3} & \overline{2} & 
\overline{1} & 1 & 2 & 3 & 4 & 5 & 6 \\ 
\overline{6} & \overline{3} & \overline{4} & \overline{5} & \overline{1} & 
\overline{2} & 2 & 1 & 5 & 4 & 3 & 6
\end{array}
\right)
\end{equation*}
and $\varepsilon (\mu )=1.$ We have $G_{\ell ,\mu }\simeq GL_{4}\times
SO_{5}.$
\end{example}

\subsection{The stabilization phenomenon\label{subsec_rq}}

We begin this paragraph with further remarks :

\bigskip

\noindent \textbf{Remarks: }

\noindent $\mathrm{(i):}$ Suppose $\varepsilon (\mu )\neq 0$. In the even
case, we have $G_{\ell ,\mu }=GL_{r_{1}}\times \cdots \times GL_{r_{p}}.$ In
the odd case and $\ell \geq p+1$, $G_{\ell ,\mu }$ is not a direct product
of linear groups since $G_{\ell ,\mu }=GL_{r_{1}}\times \cdots \times
GL_{r_{p}}\times SO_{2r_{p+1}+1}.$

\noindent $\mathrm{(ii):}$ When $\ell =2,$ we have always $\mathrm{card}%
(X^{(1)})=n=\mathrm{card}(J^{(2)}).$ Hence $\varepsilon (\mu )\neq 0$ for
all partitions $\mu .$ Observe that it does not mean that the expansion (\ref
{dec6levi}) is infinite. In fact most of the branching coefficients $[V^{%
\frak{so}_{2n+1}}(\lambda ):V^{\frak{g}_{\ell ,\mu }}(\gamma _{\ell ,\mu })]$
vanishes in this situation. Note also that we have always $G_{\ell ,\mu
}\simeq GL_{n}$ in this case.

\noindent $\mathrm{(iii):}$ We have seen that the non-zero parts of the $%
\ell $-quotient $\mu /\ell $ does not depend on the number of zero parts in $%
\mu $ (see Proposition \ref{prop_stabA}). This notably implies the stability
of the coefficients $a_{\lambda ,\mu }^{\ell ,\frak{gl}}.\;$The situation is
more subtle for the coefficients $a_{\lambda ,\mu }^{\ell ,\frak{so}}$.
Indeed, the dominant weights $\gamma _{\ell ,\mu }$ given by the previous
algorithm \textit{do} not stabilize in general when the number of zero parts
in $\mu $ increases. Let us consider for example $\mu =(1,5,5,6,9)$ and $%
\ell =2$.\ By adding parts $0$ to $\mu $, we obtain successively for the
dominant weights 
\begin{equation}
(-1,2,4,4,5),\quad (-5,-4,-4,-2,-2,1),\quad (-1,2,2,2,4,4,5)\text{ etc.}
\label{nonstable}
\end{equation}
This is not incompatible with Proposition \ref{prop-stabpleth} which assets
that $\psi _{\ell }(s_{\mu }^{\frak{so}_{2n+1}})$ stabilizes in large rank.\
In fact, this only means that, when no assumption is made on the size of $n$%
, there can exist non zero coefficients $a_{\lambda ,\mu }^{\ell ,\frak{so}%
_{2n+1}}$ in the decomposition 
\begin{equation*}
\psi _{\ell }(s_{\mu }^{\frak{so}_{2n+1}})=\sum_{\mu \in \mathcal{P}%
_{n}}a_{\lambda ,\mu }^{\ell ,\frak{so}_{2n+1}}s_{\lambda }^{\frak{so}%
_{2n+1}}
\end{equation*}
such that $a_{\lambda ,\mu }^{\ell ,\frak{so}}=0$.\ This is because the
coefficients $a_{\lambda ,\mu }^{\ell ,\frak{so}}$ coincide for $s_{\mu }^{%
\frak{so}_{2n+1}}$ and $s_{\mu }^{\frak{so}_{2n}}$ in large rank whereas $%
\psi _{\ell }(s_{\mu }^{\frak{so}_{2n+1}})\neq \psi _{\ell }(s_{\mu }^{\frak{%
so}_{2n}})$ in general when no assumption is made on the size of $n$. In the
rest of this paragraph, we are going to see that the dominant weights $\mu $
for which $\gamma _{\ell ,\mu }$ do not stabilize are such that $a_{\lambda
,\mu }^{\ell ,\frak{so}}=0,$ that is their contribution to $\psi _{\ell
}(s_{\mu }^{\frak{so}_{2n+1}})$ vanishes in large rank. Moreover, we are
going to characterize precisely these weights.

\bigskip

Suppose first $\ell =2p$ is even.\ Consider $\mu \in \mathcal{P}_{m}$ such
that $\varepsilon (\mu )\neq 0$.\ Set $\gamma _{\ell ,\mu }=(\mu
^{(1)},...,\mu ^{(p)}).\;$Write $\nu $ for the partition of $\mathcal{P}%
_{m+\ell }$ obtained by adding $\ell $ parts $0$ in $\mu $.\ For any $k\in
\{1,\ldots ,p\},$ Set $\mu ^{(k)}=(\mu _{-}^{(k)},\mu _{+}^{(k)})$ where $%
\mu _{-}^{(k)}$ (resp.\ $\mu _{+}^{(k)}$) is the sequence formed by the $%
s_{k}$ leftmost (resp.\ $r_{k}-s_{k}$ rightmost) components of $\mu ^{(k)}$
(see Section \ref{subsec-algoeven} for the notation).

\begin{lemma}
We have $\varepsilon (\nu )=\varepsilon (\mu )$. Moreover if we set $\gamma
_{\ell ,\nu }=(\nu ^{(1)},...,\nu ^{(p)})$, we obtain 
\begin{equation*}
\nu ^{(k)}=(\mu _{-}^{(k)},\alpha _{k}+1-s_{k},\alpha _{k}+1-s_{k},\mu
_{+}^{(k)})
\end{equation*}
for any $k\in \{1,\ldots ,p\},$ that is $\nu ^{(k)}$ is obtained by
inserting in $\mu ^{(k)}$ two components equal to $\alpha _{k}+1-s_{k}$.\ In
particular 
\begin{equation}
\left| \gamma _{\ell ,\nu }\right| =\left| \gamma _{\ell ,\mu }\right|
+2\sum_{k=1}^{p}\left| \alpha _{k}+1-s_{k}\right| .  \label{comparcont}
\end{equation}
\end{lemma}

\begin{proof}
Let us slightly abuse the notation by writing $I^{(k)}(\mu ),$ $J^{(k)}(\mu
),$ $I^{(k)}(\nu ),$ $J^{(k)}(\nu ),$ $k=1,\ldots ,\ell $ and $X^{(k)}(\mu ),
$ $X^{(k)}(\nu ),$ $k=1\ldots ,p$ for the sequences defined from $\mu $ and $%
\nu $ by applying the procedure of Section \ref{subsec-algoeven}.\ We define 
$\alpha _{k}(\mu )$ and $\alpha _{k}(\nu ),$ $k=1,\ldots ,p$ similarly. We
have $I^{(k)}(\nu )=\{k\}\cup (I^{(k)}(\mu )+1)$ for $k=1,\ldots ,\ell $.\
Moreover $\mathrm{card}(J^{(k)}(\nu ))=\mathrm{card}(J^{(k)}(\mu ))+2$ and $%
\alpha _{k}(\nu )=\alpha _{k}(\mu )+1$ for any $k\in \{1,\ldots ,\ell \}.$
Thus $\mathrm{card}(X^{(k)}(\nu ))=\mathrm{card}(J^{(k)}(\nu ))$ for any $%
k\in \{1,\ldots ,p\}$ and thus, $\varepsilon (\nu )=\varepsilon (\mu )$. So
it makes sense to consider $\gamma _{\ell ,\nu }=(\nu ^{(1)},...,\nu
^{(p)}).\;$It then follows by a direct application of the formulas (\ref
{mukB}) that 
\begin{equation*}
\nu ^{(k)}=(\mu _{-}^{(k)},\alpha _{k}+1-s_{k},\alpha _{k}+1-s_{k},\mu
_{+}^{(k)})
\end{equation*}
and thus $\left| \gamma _{\ell ,\nu }\right| =\left| \gamma _{\ell ,\mu
}\right| +2\sum_{k=1}^{p}\left| \alpha _{k}+1-s_{k}\right| .$
\end{proof}

\bigskip

When $\ell =2p+1$ is odd and $\gamma _{\ell ,\mu }=(\mu ^{(1)},...,\mu
^{(p)},\mu ^{(p+1})$, we can define $\nu $ similarly.\ Then, one proves that 
$\varepsilon (\nu )=\varepsilon (\mu ).\;$We have $\gamma _{\ell ,\nu }=(\nu
^{(1)},...,\nu ^{(p)},\nu ^{(p+1})$ with 
\begin{equation}
\nu ^{(k)}=(\mu _{-}^{(k)},\alpha _{k}+1-s_{k},\alpha _{k}+1-s_{k},\mu
_{+}^{(k)})\text{ for any }k=1,\ldots ,p  \label{stab_muk}
\end{equation}
and $\nu ^{(p+1)}=(0,\mu ^{(p+1)}).$ Hence (\ref{comparcont}) still holds.
With the notation of Sections \ref{subsec-algoeven} and \ref{subsec-algoodd}%
, we obtain the following stabilization theorem :

\begin{theorem}
\label{Th_LI}Consider $\mu $ a partition such that $\varepsilon (\mu )\neq 0$%
.$\;$Let $\ell $ be a positive integer. Then for any partition $\lambda $

\begin{enumerate}
\item  $a_{\lambda ,\mu }^{\ell ,\frak{so}}\neq 0$ only if $s_{k}=\alpha
_{k}+1$ for any $k=1,\ldots ,p$.

\item  In this case we have $a_{\lambda ,\mu }^{\ell ,\frak{so}}=[V^{\frak{so%
}_{2n+1}}(\lambda ):V^{\frak{g}_{\ell ,\mu }}(\gamma _{\ell ,\mu })]$ and
the non-zero components of the dominant weight $\gamma _{\ell ,\mu }$ do not
depend on the number of parts $0$ in $\lambda $ and $\mu $.
\end{enumerate}
\end{theorem}

\begin{proof}
Suppose there exists $k\in \{1,\ldots ,p\}$ such that $s_{k}\neq \alpha
_{k}+1$.\ Write $\mu (a)$ for the partition obtained by adding $a\ell $
components $0$ to $\mu .\;$By (\ref{comparcont})$,$ we have then $\left|
\gamma _{\ell ,\mu (a)}\right| \geq \left| \gamma _{\ell ,\mu }\right| +2a$%
.\ Thus, for $a$ sufficiently large, one has $\left| \gamma _{\ell ,\mu
(a)}\right| >\left| \lambda \right| $. For such $a$, we will obtain $[V^{%
\frak{so}_{2n+1}}(\lambda ):V^{\frak{g}_{\ell ,\mu }}(\gamma _{\ell ,\mu
})]=0$. Hence $[V^{\frak{so}_{2n+1}}(\lambda ):V^{\frak{g}_{\ell ,\mu
}}(\gamma _{\ell ,\mu })]$ does not coincide with a non-zero coefficient $%
a_{\lambda ,\mu }^{\ell ,\frak{so}}$. When $s_{k}=\alpha _{k}+1$ for any $%
k=1,\ldots ,p,$ the second assertion of the theorem follows from (\ref
{stab_muk}).
\end{proof}

\bigskip

\noindent \textbf{Remarks: }

\noindent $\mathrm{(i):}$ There exist very efficient procedures to compute
the branching coefficients $[V^{\frak{so}_{2n+1}}(\lambda ):V^{\frak{g}%
_{\ell ,\mu }}(\gamma _{\ell ,\mu })]$ (see \cite{Ki1}). By the previous
theorem, they permit to derive the coefficients $a_{\lambda ,\mu }^{\ell ,%
\frak{so}}$.

\noindent $\mathrm{(ii):}$ One can check that condition 1 of Theorem \ref
{Th_LI} is satisfied in Example \ref{ex-2} but fails in (\ref{nonstable})
where $s_{1}=1$ and $\alpha _{1}=2$.

\subsection{Coefficients $a_{\protect\lambda ,\protect\mu }^{\ell }$ and
restriction to Levi subgroups\label{subsec_main-TH}}

By combining the results of Sections \ref{sec-stabplethys} and \ref{sec-algo}
we derive the following theorem which expresses $a_{\lambda ,\mu }^{\ell ,%
\frak{so}}$ and $a_{\lambda ,\mu }^{\ell ,\frak{sp}}$ as branching
coefficients corresponding to restrictions to Levi subgroups.

\begin{theorem}
\label{th_a=bL}Consider $\lambda \in \mathcal{P}_{m}$ and $\ell $ a positive
integer. Let $\frak{g}$ be a symplectic or orthogonal Lie group with rank $n>%
\mathrm{max}(\ell l(\lambda ),l(\lambda ^{\prime })).$ Then we have :

\begin{enumerate}
\item  $a_{\lambda ,\mu }^{\ell ,\frak{so}}=\varepsilon (\mu )[V^{\frak{so}%
_{2n+1}}(\lambda ):V^{\frak{g}_{\ell ,\mu }}(\gamma _{\ell ,\mu })]$ where $%
\varepsilon (\mu ),$ $\frak{g}_{\ell ,\mu }$ and $\gamma _{\ell ,\mu }$ are
determined by the algorithms of Section \ref{sec-algo},

\item  $a_{\lambda ,\mu }^{\ell ,\frak{sp}}=\varepsilon (\mu ^{\prime })[V^{%
\frak{so}_{2n+1}}(\lambda ^{\prime }):V^{\frak{g}_{\ell ,\mu }}(\gamma
_{\ell ,\mu ^{\prime }})]$ if $\left| \lambda \right| $ is even and

\item  $a_{\lambda ,\mu }^{\ell ,\frak{sp}}=(-1)^{\ell -1}\varepsilon (\mu
^{\prime })[V^{\frak{so}_{2n+1}}(\lambda ^{\prime }):V^{\frak{g}_{\ell ,\mu
}}(\gamma _{\ell ,\mu ^{\prime }})]$ if $\left| \lambda \right| $ is odd.
\end{enumerate}
\end{theorem}

\begin{proof}
Assertion $1$ follows from Proposition \ref{prop-stabpleth} and Theorem \ref
{Th-decB}. By remark following Lemma \ref{lem_orthsympl_dual}, one has $%
a_{\lambda ,\mu }^{\ell ,\frak{sp}}=a_{\lambda ^{\prime },\mu ^{\prime
}}^{\ell ,\frak{so}}$ if $\left| \lambda \right| $ is even and $a_{\lambda
,\mu }^{\ell ,\frak{sp}}=(-1)^{\ell -1}a_{\lambda ^{\prime },\mu ^{\prime
}}^{\ell ,\frak{so}}$ otherwise which proves assertion $2$. Note that the
assumption $n>\mathrm{max}(\ell l(\lambda ),l(\lambda ^{\prime }))$ suffices
to guarantee that $\lambda ^{\prime }$ belongs to $\mathcal{P}_{n}$.
\end{proof}

\bigskip 

\noindent In the sequel, we will assume for simplicity $\ell \geq 2$ and $%
n\geq \ell \left| \lambda \right| $ which implies the condition $n>\mathrm{%
max}(\ell l(\lambda ),l(\lambda ^{\prime }))$.

\section{Splitting $V^{\frak{g}}(\protect\lambda )^{\otimes 2}$ into its
symmetric and antisymmetric parts\label{sec_spli}}

\subsection{Decomposition of the plethysms $p_{2}\circ s_{\protect\lambda }^{%
\frak{g}}$}

Consider a partition $\lambda \in \mathcal{P}_{m}.\;$According to Theorem 
\ref{th_a=bL}, we have with the notation of Sections \ref{subsec_ple_shur}
and \ref{subsec-stab} 
\begin{eqnarray}
p_{2}\circ s_{\lambda }^{\frak{gl}_{n}} &=&\sum_{\mu \in \mathcal{P}%
_{n}}\varepsilon (\mu )c_{(\mu ^{(0)},\mu ^{(1)})}^{\lambda }s_{\mu }^{\frak{%
gl}_{n}},  \label{ide} \\
p_{2}\circ s_{\lambda }^{\frak{so}} &=&\sum_{\mu \in \mathcal{P}%
_{n}}\varepsilon (\mu )[V^{\frak{so}_{2n+1}}(\lambda ):V^{\frak{gl}%
_{n}}(\gamma _{\mu })]s_{\mu }^{\frak{so}}\text{ and }p_{2}\circ s_{\lambda
}^{\frak{sp}}=(-1)^{\left| \lambda \right| }p_{2}\circ s_{\lambda ^{\prime
}}^{\frak{so}}  \notag
\end{eqnarray}
for any $n\geq 2\left| \lambda \right| .\;$Here we have written for short $%
\gamma _{\mu }$ for $\gamma _{2,\mu }$ and $V^{\frak{gl}_{n}}(\gamma _{\mu
}) $ instead of $V^{\frak{L}_{2,\mu }}(\gamma _{\mu })$ (see Remark $\mathrm{%
(ii)}$ of Section \ref{subsec_rq}). Since $n\geq m$ and $\gamma _{\mu
}=(\gamma ^{-},\gamma ^{+})$ belongs to $\widetilde{\mathcal{P}}_{n}$ we
have the following decomposition (see \cite{Ki1}) : 
\begin{equation}
\lbrack V^{\frak{so}_{2n+1}}(\lambda ):V^{\frak{gl}_{n}}(\gamma _{\mu
})]=\sum_{\delta ,\xi \in \mathcal{P}_{n}}c_{\delta ,\xi }^{\lambda
}c_{\gamma ^{-},\gamma ^{+}}^{\xi }.  \label{branB_lr}
\end{equation}

\subsection{Symmetric and antisymmetric parts of $V^{\frak{g}}(\protect%
\lambda )^{\otimes 2}$}

\noindent Consider $\lambda \in \mathcal{P}_{m}.\;$By Propositions \ref
{stab_pleth_schur} and \ref{prop-stabpleth}, for any rank $n\geq 2\left|
\lambda \right| $ the plethysms $h_{2}\circ s_{\lambda }^{\frak{g}}$ and $%
e_{2}\circ s_{\lambda }^{\frak{g}}$ stabilize. Set 
\begin{equation*}
h_{2}\circ s_{\lambda }^{\frak{g}}=\sum_{\mu \in \mathcal{P}_{n}}m_{\lambda
,\mu }^{\frak{G},+}s_{\mu }^{\frak{g}}\text{\qquad and\qquad }e_{2}\circ
s_{\lambda }^{\frak{g}}=\sum_{\mu \in \mathcal{P}_{n}}m_{\lambda ,\mu }^{%
\frak{G,-}}s_{\mu }^{\frak{g}}
\end{equation*}
where $\frak{G=gl,so}$ or $\frak{sp}$ respectively when $\frak{g=gl}_{n},%
\frak{so}_{N}$ or $\frak{sp}_{2n}.\;$Recall that $h_{2}\circ s_{\lambda }^{%
\frak{g}}$ and $e_{2}\circ s_{\lambda }^{\frak{g}}$ are the characters of $%
S^{2}(V^{\frak{g}}(\lambda )$ and $\Lambda ^{2}(V^{\frak{g}}(\lambda ).$ By
using (\ref{relfund}) and Theorem \ref{th_a=bL}, we obtain for any rank $%
n\geq 2\left| \lambda \right| $%
\begin{eqnarray*}
m_{\lambda ,\mu }^{\frak{gl},\pm } &=&\frac{1}{2}(c_{\lambda ,\lambda }^{\mu
},\pm \varepsilon (\mu )c_{(\mu ^{(0)},\mu ^{(1)})}^{\lambda }),\text{ } \\
m_{\lambda ,\mu }^{\frak{so},\pm } &=&\frac{1}{2}(d_{\lambda ,\lambda }^{\mu
},\pm \varepsilon (\mu )[V^{\frak{so}_{2n+1}}(\lambda ):V^{\frak{gl}%
_{n}}(\gamma _{\mu })],\text{ } \\
m_{\lambda ,\mu }^{\frak{sp},\pm } &=&\frac{1}{2}(d_{\lambda ^{\prime
},\lambda ^{\prime }}^{\mu ^{\prime }},\pm (-1)^{\left| \lambda \right|
}\varepsilon (\mu ^{\prime })[V^{\frak{so}_{2n+1}}(\lambda ^{\prime }):V^{%
\frak{gl}_{n}}(\gamma _{\mu ^{\prime }})]
\end{eqnarray*}
where the coefficients $d_{\lambda ,\lambda }^{\mu }$ are the multiplicities
appearing in Proposition \ref{prop_inde}. Now these multiplicities can be
expressed in terms of the Littlewood coefficients \cite{Ki}.\ Namely we have 
$d_{\lambda ,\lambda }^{\mu }=\sum_{\delta ,\xi ,\eta }c_{\delta ,\xi }^{\mu
}c_{\delta ,\eta }^{\lambda }c_{\xi ,\eta }^{\lambda }.$ In particular we
recover the equality $d_{\lambda ^{\prime },\lambda ^{\prime }}^{\mu
^{\prime }}=$ $d_{\lambda ,\lambda }^{\mu }$ since $c_{\delta ,\eta
}^{\gamma }=c_{\delta ^{\prime },\eta ^{\prime }}^{\gamma ^{\prime }}$ for
any partitions $\delta ,\eta $ and $\gamma .$ By using (\ref{branB_lr}),
this thus permits to express the multiplicities appearing in the symmetric
and antisymmetric parts of $V^{\frak{g}}(\lambda )^{\otimes 2}$ in terms of
the Littlewood-Richardson coefficients. Note that formulas for computing the
plethysms $h_{2}\circ s_{\lambda }^{\frak{g}}$ and $e_{2}\circ s_{\lambda }^{%
\frak{g}}$ were introduced without a complete proof by Littlewood in \cite
{Li2}.

\begin{proposition}
\label{prop-symanti}With the above notation we have for any rank $n\geq
2\left| \lambda \right| $%
\begin{eqnarray*}
m_{\lambda ,\mu }^{\frak{gl},\pm }=\frac{1}{2}(c_{\lambda ,\lambda }^{\mu
}\pm \varepsilon (\mu )c_{(\mu ^{(0)},\mu ^{(1)})}^{\lambda }), \\
m_{\lambda ,\mu }^{\frak{so,}\pm }=\frac{1}{2}\left( \sum_{\delta ,\xi ,\eta
\in \mathcal{P}_{n}}c_{\delta ,\xi }^{\mu }c_{\delta ,\eta }^{\lambda
}c_{\xi ,\eta }^{\lambda }\pm \varepsilon (\mu )\sum_{\delta ,\xi \in 
\mathcal{P}_{n}}c_{\delta ,\xi }^{\lambda }c_{\gamma ^{-},\gamma ^{+}}^{\xi
}\right) \\
m_{\lambda ,\mu }^{\frak{sp},\pm }=\frac{1}{2}\left( \sum_{\delta ,\xi ,\eta
\in \mathcal{P}_{n}}c_{\delta ,\xi }^{\mu }c_{\delta ,\eta }^{\lambda
}c_{\xi ,\eta }^{\lambda }\pm (-1)^{\left| \lambda \right| }\varepsilon (\mu
^{\prime })\sum_{\delta ,\xi \in \mathcal{P}_{n}}c_{\delta ,\xi }^{\lambda
}c_{\kappa ^{-},\kappa ^{+}}^{\xi }\right)
\end{eqnarray*}
where $\gamma _{\mu }=(\gamma ^{-},\gamma ^{+})$ and $\gamma _{\mu ^{\prime
}}=(\kappa ^{-},\kappa ^{+}).$
\end{proposition}

\end{document}